\pdfoutput=1 
\documentclass[12pt]{article}

\usepackage[margin=1in]{geometry}  
\usepackage{graphicx}              
\usepackage{amsmath}               
\usepackage{amsfonts}              
\usepackage{amsthm}                
\usepackage{url}                   

\DeclareMathSizes{12}{12}{11}{11}     

\newcommand{\sign}[1]{\textsf{#1}} 

\title{\textbf{Non-linear mass-spring system for large soft tissue deformations modeling}} 

\author{
\large
\textsc{Nikolaev S. N.}\\[2mm] 
\normalsize Saint-Petersburg State University \\ 
\normalsize mailto:  ser.niev@rambler.ru 
\vspace{-5mm}
}
\date{}

\begin{document}

article information: author -- Nikolaev S.N.; title -- Non-linear mass-spring system for large soft tissue deformations modeling;
journal title -- Scientific and technical journal of information technologies, mechanics and optics; date -- September-October, 2013;
number -- 5(87); pages -- 88-94; url -- \url{http://ntv.ifmo.ru/en/article/4322/article_4322.htm}.

\maketitle

\begin{abstract}
  Implant placement under soft tissues operation is described. In this operation tissues can reach such deformations that nonlinear properties are appeared. A mass-spring model modification for modeling nonlinear tissue operation is developed. A method for creating elasticity module using splines is described. For Poisson ratio different stiffness for different types of springs in cubic grid is used. For stiffness finding an equation system that described material tension is solved. The model is verified with quadratic sample tension experiment. These tests show that sample tension under external forces is equal to defined nonlinear elasticity module. The accuracy of Poisson ratio modeling is thirty five percent that is better the results of available ratio modeling method.

\end{abstract}

\section{Introduction}

\hspace{6mm}Nowadays medicine develops rapidly. Surgeons can perform more quality and safer operations with modern equipment, materials and processes. One of the most popular operations among all is operation where implant is placed under soft tissues. Such operations are usually performed in plastic and reconstructive surgery. The aims of plastic surgery are improving the appearance of patient and body correction with right proportions. In reconstructive cases the implantation is performed for donor tissues growing. To stimulate growth the expanders are usually used. Nowadays surgeons perform such operations using only their experience. So an ability to visualize operations results could help surgeons to select right implant or estimate donor tissue sizes.

Finite element method is the most complete and qualitative method to solve biomechanical tasks. However surgeons need to watch modeling result as fast as possible so executable speed is also a key point. Because of complex calculations finite element model works slowly. So nowadays together with it the new mass-spring model \cite{kass}, the main criteria of which is modeling speed, also becomes widely used. This approach was successfully implemented to model different biomechanical tasks: heart \cite{jarrouse, wang}, lungs \cite{wang} and muscles \cite{nedel}.

To model the described above operation an implant is placed under soft tissues and stretched with them \cite{nikolaev}. In this case soft tissues undergo tensile deformations. To calculate the maximal deformation sizes we can set a hemisphere as upper implant restriction (because every implant height is less than its radius) and plane as lower stretchable tissues restriction. Then the area of stretchable surface is equal to area of a circle:
$$
S_1 = \pi * R
$$
where R equals implant radius. Surface area after extension is equal to area of a hemisphere:
$$
S_2 = 2 * \pi * R
$$
Thus tissues relative tension ratio can achieve a value of 1.0. It is known that in such deformations soft tissues have nonlinear elasticity that cannot be neglected \cite{fung}.

In common mass-spring model springs have linear stress-strain relations \cite{kass}. Several works describe this model improvements where spring relations are considered as second and third order polynomials \cite{jarrouse} or as a piecewise linear function \cite{basafa}. However nonlinear tissues behavior curves are much more complicated than described polynomials. And piecewise linear functions approximate exact these curves only in sections borders. Also when a tissue is expand in one direction it tends to compress in the other two directions perpendicular to the direction of expansion. This is described by Poisson ratio. Poisson ratio is modeled with genetic algorithm \cite{bianchi} or analytical approach \cite{baudet}. For genetic algorithm a different measurements of artificial samples are need, which is often impossible to get. In analytical approach a shear module is used for diagonal springs, which has right behavior only for small deformations.

As a result we need a new method for modeling nonlinear elasticity and Poisson ratio. The aim of this paper is to develop mass-spring model modification for modeling tissues deformation in diapasons of relative tension ratio between 0.0 and 1.0 values.

\section{Nonlinearity modeling}

\hspace{6mm}In mass-spring model an object is considered as points set that describes object mass and space position. The points in set are connected with springs (fig. 1)

\renewcommand{\figurename}{Fig.}
\begin{figure}[ht!]
\centering
\includegraphics[width=50mm]{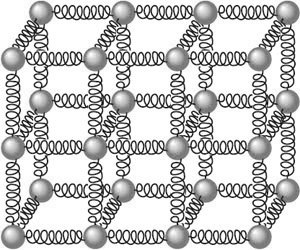}
\caption{\sign{Mass-spring model}}
\label{overflow}
\end{figure}

Spring action formula is declared as:
$$
f = k * (|x_{12}| - l_{12}) * \frac{x_{12}}{|x_{12}|}
$$
here $f$ -- force acted on the one of vertexes, which is connected to spring (the force acted on the second vertex has the same magnitude but opposite direction), $x_{12}$ -- current distance between vertexes, $l_{12}$ -- initial spring length, $k$ -- spring stiffness coefficient. We can also introduce another form for spring formula:
$$
f = Sq * E * \frac{\Delta L}{L}
$$
where $E$ -- Young module, $Sq$ -- topological element size, $L$  -- initial spring length, $\Delta L$ -- absolute spring elongation magnitude. This formula is more preferred because it is valid for spring with any size. In nonlinear case this equation is transferred to more common expression:
$$
f = Sq * Ef(\frac{\Delta L}{L})
$$
Function $Ef$ describes relation between stress and strain.

To construct $Ef$ function we need to mark some points on initial relation chart. The points are connected with segment lines or splines. As a result we get a piecewise smooth function that consists of polynomial, trigonometric or exponential parts.

\section{Poisson ratio modeling}

\hspace{6mm}This work covers only modeling in two-dimensional space. For Poisson ratio modeling we take square grid. Grid springs are divided into edge springs (located on square edges) and diagonal springs (located on square diagonals). Consider such square with the length size $L_{\text{edge}}$ (fig. 2). Suppose it elongates along some axe direction on a value of $\Delta L_{\text{edge}}$.

\begin{figure}[ht!]
\centering
\includegraphics[width=90mm]{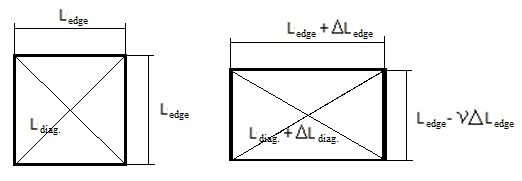}
\caption{\sign{Square grid elongation}}
\label{overflow}
\end{figure}

In transversal directions square compresses on a next value:
$$
\Delta L_{\text{comp.}} = \epsilon_{\text{comp.}} * L_{\text{edge}} = -\nu * \epsilon_{\text{elong.}} * L_{\text{edge}} = -\nu * \Delta L_{\text{edge}}
$$
here $\nu$  -- Poisson ratio. In paper \cite{baudet} authors express relative springs coefficients in terms of shear and elongation modules. However we need to restrict deformations only with elongations because springs do not have shear properties. So the constructed equations system looks like this:
\begin{equation}
\begin{cases}
\frac{F}{L_{\text{edge}}} = Ef(\frac{\Delta L_{\text{edge}}}{L_{\text{edge}}}) = k_1 * Ef(\frac{\Delta L_{\text{edge}}}{L_{\text{edge}}}) + d_1 * Ef(\frac{\Delta L_{\text{diag.}}}{L_{\text{diag.}}}) * \frac{L_{\text{edge}} + \Delta L_{\text{edge}}}{L_{\text{diag.}} + \Delta L_{\text{diag.}}} \\
k_2 * Ef(\frac{-\nu * \Delta L_{\text{edge}}}{L_{\text{edge}}}) + d_1 * Ef(\frac{\Delta L_{\text{diag.}}}{L_{\text{diag.}}}) * \frac{L_{\text{edge}} - \nu * \Delta L_{\text{edge}}}{L_{\text{diag.}} + \Delta L_{\text{diag.}}} = 0
\end{cases}
\end{equation}
where $k_1$ -- is the unknown edge springs coefficient, $d_1$ -- the unknown diagonal springs coefficient, $k_2$ describes the case, when edge spring is compressed.

Express $d_1$ from the second equation of system (1) and substitute it in the first equation. As a result we get:
\begin{equation}
Ef(\frac{\Delta L_{\text{edge}}}{L_{\text{edge}}}) = k_1 * Ef(\frac{\Delta L_{\text{edge}}}{L_{\text{edge}}}) - k_2 * Ef(\frac{-\nu * \Delta L_{\text{edge}}}{L_{\text{edge}}}) * \frac{L_{\text{edge}} + \Delta L_{\text{edge}}}{L_{\text{edge}} - \nu * \Delta L_{\text{edge}}}
\end{equation}
In equation (2) there are two unknown coefficients, so one of them can have any value. In this work we set to $k_2$ the value 1.0 and express $k_1$ via $k_2$ in equation. Thus the common formula for forces that define edge springs action on connected vertices is look like this:
\begin{equation}
f_{\text{edge}} =
\begin{cases}
Sq * (Ef(\frac{\Delta L_{\text{edge}}}{L_{\text{edge}}}) + \frac{ Ef(\frac{-\nu * \Delta L_{\text{edge}}}{L_{\text{edge}}})}{\frac{(1 + \nu) * L_{\text{edge}}}{L_{\text{edge}} + \Delta L_{\text{edge}}} - \nu } ) & ,\Delta L_{\text{edge}} \ge 0 \\
Sq * Ef(\frac{\Delta L_{\text{edge}}}{L_{\text{edge}}}) & ,\Delta L_{\text{edge}} < 0
\end{cases}
\end{equation}
$\Delta L_{\text{diag.}}$ and $L_{\text{diag.}}$ describe diagonal springs lengths in system (1). Using Pythagorean Theorem conditions we can express them via edge springs lengths:
\begin{equation}
L_{\text{diag.}} = L_{\text{edge}} * \sqrt{2}
\end{equation}
\begin{equation}
(L_{\text{diag.}} + \Delta L_{\text{diag.}}) = \sqrt{(L_{\text{edge}} + \Delta L_{\text{edge}})^2 + (L_{\text{edge}} - \nu * \Delta L_{\text{edge}})^2}
\end{equation}
$\Delta L_{\text{edge}}$ can be found via known diagonal springs lengths from equations (4), (5):
$$
\Delta L_{\text{edge}} = \frac{L_{\text{diag.}}}{\sqrt{2}} * J = L_{\text{edge}} * J 
$$
where $J$ is defined as:
$$
J = \frac{\nu - 1 + \sqrt{2 * (1 + \nu^2) * \left(\frac{L_{\text{diag.}} + \Delta L_{\text{diag.}}}{L_{\text{diag.}}} \right)^2 - (1 + 2\nu + \nu^2)}}{(1 + \nu^2)}
$$
Substituting $\Delta L_{\text{edge}}$ in the second equation of system (1) we get the next formula ( $k_2$ is already equal to 1.0):
$$
d_1 * Ef(\frac{\Delta L_{\text{diag.}}}{L_{\text{diag.}}}) = \frac{-\sqrt{2} * Ef(-\nu * J) * (L_{\text{diag.}} + \Delta L_{\text{diag.}})}{(1 - \nu * J) * L_{\text{diag.}}}
$$
Set of equations (1) does not contain diagonal springs compression behavior. So for diagonal springs compression we use the same formula as for edges springs compression. As a result we can write the common formula for forces that define diagonal springs action on connected vertices:
\begin{equation}
f_{\text{diag.}} =
\begin{cases}
Sq * \frac{-\sqrt{2} * Ef(-\nu * J) * (L_{\text{diag.}} + \Delta L_{\text{diag.}})}{(1 - \nu * J) * L_{\text{diag.}}} & ,\Delta L_{\text{diag.}} \ge 0 \\
Sq * Ef(\frac{\Delta L_{\text{diag.}}}{L_{\text{diag.}}}) & ,\Delta L_{\text{diag.}} < 0
\end{cases}
\end{equation}
The derived formulas are substituted instead of linear spring formulas. In conclusion we should mark that function $Ef$ , which describe nonlinear elasticity, and Poisson ratio value are only necessary for spring forces calculation in formulas (3) and (6) in new model.

\section{Elongation deformation modeling results}

\hspace{6mm}We performed some experiments after new model constructing. During experiments models with different tissues materials were elongated. We took squared sample with edge size equal to 4 unit of lengths, aligned it along axes and divided one into grid with edge size equal to 1 unit of length.

This model was primarily developed to simulate soft tissues deformations above implant, so we use elasticity modulus of soft tissues as source data. In the simplest case implant is placed under two soft tissues: skin and adipose. Thus we use these tissues data for experiments. As a separate experiment we calculate stress values for different strains with given $Ef$ function. The results are presented on charts 1 and 3 (skin '$Ef$' and adipose '$Ef$').

We perform elongation experiments exactly as the description below. The two elongation forces were acted on sample along some axe in opposite directions. These forces were applied to each sample bound point in given direction. The values of forces were calculated such that forces with the same magnitude acted on each two vertexes of each grid element. In other words the doubled forces were acted on vertexes that owned to different squared grid elements. In this case the average stress which elongates the sample is equal to the half of given force value. After 15000 iterations we use bound box to calculate sizes along axes of the stretched sample. Using initial and result sizes values we compute the relative strains in elongation and transversal directions.

The experiment was repeated several times. For each considered material we began deformations with small forces. Then in each next experiment we increase forces values. The experiments were repeated until relative strain in elongation was less than 1.0.

\renewcommand{\figurename}{Chart}
\setcounter{figure}{0}
\begin{figure}[ht!]
\centering
\includegraphics[width=130mm]{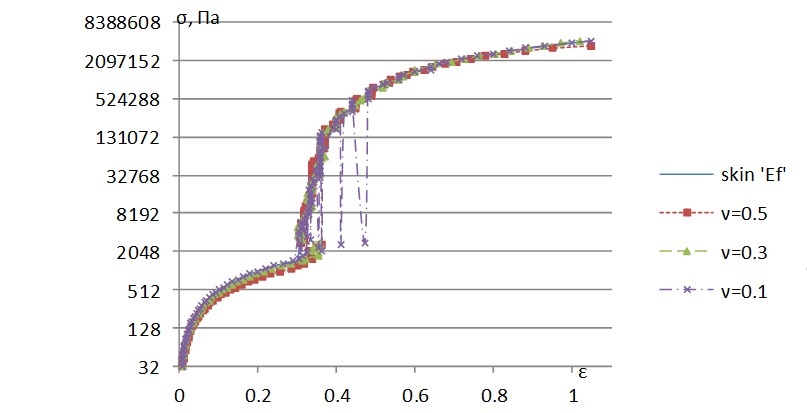}
\caption{\sign{Nonlinear elasticity skin modulus experimental results for given Poisson ratio values. skin '$Ef$' -- nonlinear elasticity modulus values}}
\label{overflow}
\end{figure}
\smallskip
\begin{figure}[ht!]
\centering
\includegraphics[width=130mm]{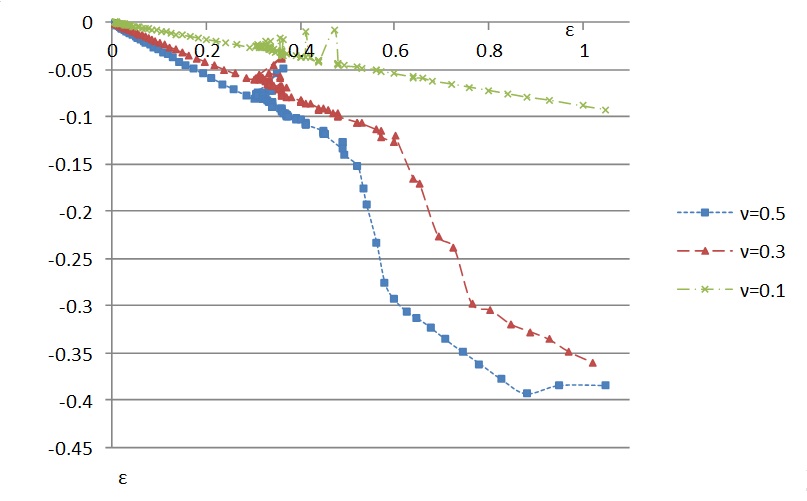}
\caption{\sign{Poisson ratio experimental results for skin}}
\label{overflow}
\end{figure}

Firstly we perform experiments with skin. Its nonlinear elasticity can be split onto three sections in relative elongation diapason between 0 and 1 \cite{hendriks}. It is constant during small deformations. Further elasticity begins to increase after strain value 0.4. When strain value becomes nearly 0.7 elasticity ends increasing and become constant again. Thus we split $Ef$ function onto three parts. First and third parts were approximated with first order functions. Second part was defined as spline that connects linear chunks. Chart 1 shows $Ef$ values and experimental relations between stress and strain for given Poisson ratio values (which are marked in legend). Chart 2 shows experimental relation between compress and elongation for given Poisson ratio values.

\begin{figure}[ht!]
\centering
\includegraphics[width=130mm]{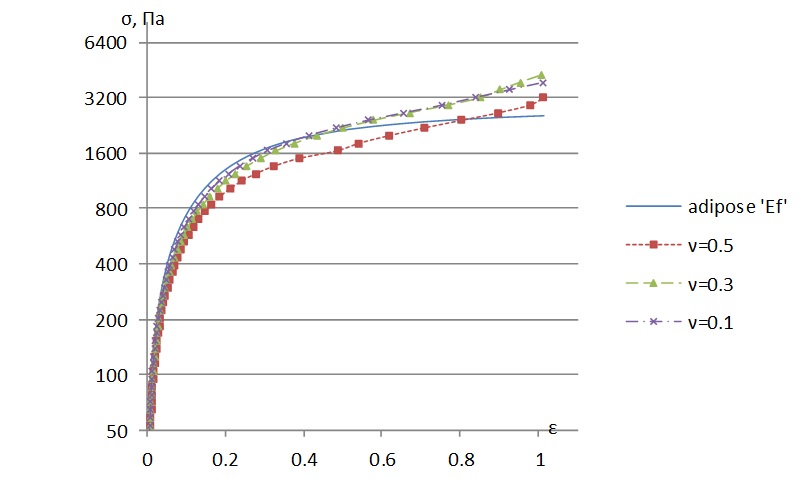}
\caption{\sign{Nonlinear elasticity adipose modulus experimental results for given Poisson ratio values. adipose '$Ef$' -- nonlinear elasticity modulus values}}
\label{overflow}
\end{figure}
\smallskip
\begin{figure}[ht!]
\centering
\includegraphics[width=130mm]{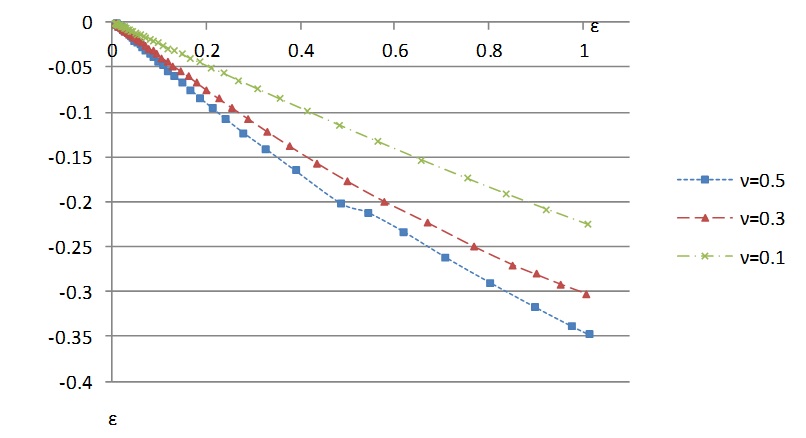}
\caption{\sign{Poisson ratio experimental results for adipose}}
\label{overflow}
\end{figure}

Adipose properties are given in \cite{geerlings}. Its elasticity is constant during small deformations. Then it decreases slowly. And when strain value becomes more 0.6 adipose elasticity stress becomes constant and does not depend anymore on material elongation. Scaled arctangent function approximate adipose elasticity behavior the most accurately. Chart 3 shows nonlinear elasticity modulus values (adipose $Ef$) and experimental relations between stress and strain for given Poisson ratio values. Chart 4 shows experimental relation between compress and elongation for adipose material.

Charts 1 and 3 demonstrate us that experiments results coincide with specified nonlinear elasticity modulus values. Some mismatching is appeared on chart 3 when strain is large. Also chart 1 shows us large error when Poisson ratio equals 0.1 and strain approximately equals 0.5. But in general experimental nonlinear modulus corresponds to the given function $Ef$.

Referring to charts 2 and 4 we can see that relation between compression and elongation is less the given Poisson ratio. For larger Poisson ratios approximation error also becomes larger. But in general relative accuracy nowhere exceeds thirty five percent.

In conclusion we want to say that experimental results for proposed method repeat forms of nonlinear elasticity curves in contrast with approaches form \cite{jarrouse, basafa}. Also Poisson ratio modeling is more accurate than in work \cite{baudet}. Elasticity of most soft tissues is similar to one of skin and fat. Moreover Poisson ratio of soft tissues usually has value between 0.35 and 0.5. So the given approach rather exact approximates this soft tissues elongation behavior. Thus this method can be used for modeling tissues under which different implants or expanders are placed.

\section{Conclusion}

\hspace{6mm}In this work we created the mass-spring model modification. The new model allows simulating material with nonlinear elasticity modulus and Poisson ratio for relative strains between 0 and 1. Also we did tests for two elasticity curves and different Poisson ratio values. These tests gave us satisfactory results.

The obtained model can simulate tissues behavior in operation where implant is placed under soft tissues \cite{nikolaev}. Also in general this method can describe any material with nonlinear elasticity modulus in given relative elongation diapason.

In future work we are going to develop model modifications for three-dimensional space.

\end{document}